\documentclass[psamsfonts]{amsart}

\usepackage{amssymb,amsfonts}
\usepackage[all,arc]{xy}
\usepackage{enumerate}
\usepackage{mathrsfs}
\usepackage{verbatim}

\newcommand{\C}{\mathbb{C}}

\newcommand{\N}{\mathbb{N}}

\newcommand{\R}{\mathbb{R}}

\newcommand{\sgn}{\mathrm{sgn}}

\newcommand{\T}{\mathbb{T}}

\newcommand{\Z}{\mathbb{Z}}

\newtheorem{thm}{Theorem}[section]

\newtheorem{fact}[thm]{Fact}

\theoremstyle{definition}
\newtheorem{defn}[thm]{Definition}

\theoremstyle{remark}

\makeatletter
\let\c@equation\c@thm
\makeatother
\numberwithin{equation}{section}

\title[The Plancherel Formula]{The Plancherel Formula, the Plancherel
  Theorem, and the Fourier Transform of Orbital Integrals}

\author[Herb--Sally]{Rebecca A. Herb (University of Maryland) and Paul
  J. Sally, Jr. (University of Chicago)}

\date{\today}

\begin{document}

\begin{abstract}
We discuss various forms of the Plancherel Formula and the Plancherel
Theorem on reductive groups over local fields.
\end{abstract}

\maketitle

\begin{center}
{\it Dedicated to Gregg Zuckerman on his 60th birthday}
\end{center}

\ \\

%% get rid of in-text paper titles in some places?

\section{Introduction}

The classical Plancherel Theorem proved in 1910 by Michel Plancherel
can be stated as follows:
\begin{thm}
Let $f\in L^2(\R )$ and define $\phi _n:\R \rightarrow \C$ for $n\in
\N$ by
\[
\phi _n(y)=\frac{1}{\sqrt{2\pi}}\int _{-n}^nf(x)e^{iyx}dx.
\]
The sequence $\phi _n$ is Cauchy in $L^2(\R )$ and we write $\phi
=\lim _{n\to \infty}\phi _n\text{ (in }L^2\text{)}$.  Define $\psi
_n:\R \rightarrow \C$ for $n\in \N$ by
\[
\psi _n(x)=\frac{1}{\sqrt{2\pi}}\int _{-n}^n\phi (y)e^{-iyx}dy.
\]
The sequence $\psi _n$ is Cauchy in $L^2(\R )$ and we write $\psi
=\lim _{n\to \infty}\psi _n\text{ (in }L^2\text{)}$.  Then,
\[
\psi =f\text{ almost everywhere, and }\int _\R \left| f(x)\right| ^2dx=\int _\R \left| \phi (y)\right| ^2dy.
\] 
\end{thm}
This theorem is true in various forms for any locally compact abelian group.  It is often proved by starting with $f\in L^1(\R )\cap L^2(\R )$, but it is really a theorem about square integrable functions.

There is also a ``smooth'' version of Fourier analysis on $\R$,
motivated by the work of Laurent Schwartz, that leads to the
Plancherel Theorem.
\begin{defn}[The Schwartz Space]
The \emph{Schwartz space}, $\mathcal{S}(\R )$, is the collection of complex-valued functions $f$ on $\R$ satisfying:
\begin{enumerate}
\item $f\in C^\infty (\R )$.
\item $f$ and  all its derivatives vanish at infinity faster
  than any polynomial.  That is, $\lim _{|x|\to
    \infty}|x|^kf^{(m)}(x)=0$ for all $k,m\in \N$.
\end{enumerate}
\end{defn}
\begin{fact}
The Schwartz space has the following properties:
\begin{enumerate}
\item The space $\mathcal{S}(\R )$ is dense in $L^p(\R )$ for $1\leq p<\infty$.
\item The space $\mathcal{S}(\R )$ is not dense in $L^\infty (\R )$.
\item The space $\mathcal{S}(\R )$ is a vector space over $\C$.
\item The space $\mathcal{S}(\R )$ is an algebra under both pointwise multiplication and convolution.
\item The space $\mathcal{S}(\R )$ is invariant under translation.
\end{enumerate}
\end{fact}

For $f\in \mathcal{S}(\R )$, we define the Fourier transform as usual by
\[
\widehat{f}(y)=\frac{1}{\sqrt{2\pi}}\int _\R f(x)e^{iyx}dx.
\]
Of course, there are no convergence problems here, and we have
\[
f(x)=\frac{1}{\sqrt{2\pi}}\int _\R \widehat{f}(y)e^{-iyx}dy.
\]
This leads to the Plancherel Theorem for functions in $\mathcal{S}(\R
)$ by setting $\widetilde{f}(x)=\overline{f(-x)}$ and considering
$f\ast \widetilde{f}$ at $0$.  Using the fact that the Fourier
transform carries convolution product to function product, we have
\[
\left\| f\right\| ^2=\left[ f\ast \widetilde{f}\right]
(0)=\frac{1}{\sqrt{2\pi}}\int _\R \widehat{f\ast
  \widetilde{f}}(y)dy=\left\| \widehat{f}\right\| ^2.
\]

It is often simpler to work on the space $C_c^\infty (\R )$ of
complex-valued, compactly supported, infinitely differentiable
functions on $\R$.  However, nonzero functions in $C_c^\infty (\R )$
do not have Fourier transforms in $C_c^\infty (\R )$.  On the other
hand, the Fourier transform is an isometric isomorphism from
$\mathcal{S}(\R )$ to $\mathcal{S}(\R )$.

The spaces $C_c^\infty (\R )$ and $\mathcal{S}(\R )$ can be turned into topological vector spaces so that the embedding from $C_c^\infty (\R )$ into $\mathcal{S}(\R )$ is continuous.  However, the topology on $C_c^\infty (\R )$ is not the relative topology from $\mathcal{S}(\R )$.  A continuous linear functional on $C_c^\infty (\R )$ is a \emph{distribution} on $\R$, and this distribution is \emph{tempered} if it can be extended to a continuous linear functional on $\mathcal{S}(\R )$ with the appropriate topology.  This situation will arise again in our discussion of the Plancherel Formula on reductive groups.

Work on the Plancherel Formula for non-abelian groups began in earnest
in the late 1940s.  There were two distinct approaches.  The first,
for separable, locally compact, unimodular groups, was pursued by
Mautner \cite{mautner}, Segal \cite{segal}, and others.  The second,
for semisimple Lie groups, was followed by Gel$'$fand--Naimark
\cite{gelfandnaimark}, and Harish--Chandra \cite{hcplanchcss}, along
with others.  Segal's paper \cite{segal} and Mautner's paper
\cite{mautner} led eventually to the following statement (see
\cite{folland}, Theorem 7.44).
\begin{thm}
Let $G$ be a separable, unimodular, type I group, and let $dx$ be a fixed
Haar measure on $G$.  There exists a positive measure $\mu$ on
$\widehat{G}$ (determined uniquely up to a constant that depends only
on $dx$) such that, for $f\in L^1(G)\cap L^2(G)$, $\pi (f)$ is a
Hilbert--Schmidt operator for $\mu$-almost all $\pi \in \widehat{G}$,
and
\[
\int _G\left| f(x)\right| ^2dx=\int _{\widehat{G}}\left\| \pi (f)\right\| _{\text{HS}}^2d\mu (\pi ).
\]
\end{thm}
Here, of course, $\widehat{G}$ denotes the set of equivalence classes
of irreducible unitary representations of $G$.

At about the same time, Harish-Chandra stated the following theorem in his paper \emph{Plancherel Formula for Complex Semisimple Lie Groups}.
\begin{thm}
Let $G$ be a connected, complex, semisimple Lie group.  Then,
for $f\in C_c^\infty(G)$,
\[
f(1)=\lim _{H \to 0}\prod _{\alpha \in P}
D_\alpha\overline{D_\alpha}\left[ e^{\rho (H)+\overline{\rho
      (H)}}\int _{K\times N}f\left( u\exp
  (H)nu^{-1}\right) dudn\right] .
\]
\end{thm}
An explanation of the notation here can be found in
\cite{hcplanchcss}.  We do note two things.  First of all, $f$ is
taken to be in $C_c^\infty (G)$, and the formula for $f(1)$ is the
limit of a differential operator applied to what may be regarded as a
Fourier inversion formula for the orbital integral over a conjugacy
class of $\exp(H)$ in $G$.  It should also be mentioned that
not all irreducible unitary representations are contained in the
support of the Plancherel measure for complex semisimple Lie groups.
In particular, the complementary series are omitted.

In this note, we will trace the evolution of the Plancherel Formula
over the past sixty years.  For real groups, we observe that the
original Plancherel Formula and the Fourier inversion formula
ultimately became a decomposition of the Schwartz space into
orthogonal components indexed by conjugacy classes of Cartan
subgroups.  While this distinction might not have been clear for real
semisimple Lie groups, it certainly appeared in the development of the
Plancherel Theorem for reductive $p$-adic groups by Harish-Chandra in
his paper \emph{The Plancherel Formula for Reductive $p$-adic Groups}
in \cite{hccollected4}.  See also the papers of Waldspurger
\cite{waldsplanch} and Silberger \cite{silbergerplanch},
\cite{silbergerplancherrata}.  For $p$-adic groups, the lack of
information about irreducible characters and suitable techniques for
Fourier inversion has made the derivation of an explicit Plancherel
Formula very difficult.

In this paper, the authors have drawn extensively on the perceptive
description of Harish-Chandra's work by R. Howe, V. S. Varadarajan,
and N. Wallach (see \cite{hccollected}).  The authors would like to
thank Jonathan Gleason and Nick Ramsey for their assistance in
preparing this paper.  We also thank David Vogan for his valuable
comments on the first draft.

\section{Orbital Integrals and the Plancherel Formula}

Let $G$ be a reductive group over a local field.  For $\gamma \in G$,
let $G_\gamma$ be the centralizer of $\gamma$ in $G$.  Assume
$G_\gamma$ is unimodular.  For $f$ ``smooth'' on $G$, define
\[
\Lambda _\gamma (f)=\int _{G/G_\gamma}f\left( x\gamma x^{-1}\right) d\dot{x},
\]
with $d\dot{x}$ a $G$-invariant measure on $G/G_\gamma$.

Then, $\Lambda _\gamma$ is an invariant distribution on $G$,
that is, $\Lambda _\gamma (f)=\Lambda _\gamma (^yf)$ where $^yf(x)=f\left(
yxy^{-1}\right)$ for $y\in G$.  A major problem in harmonic analysis
on reductive groups is to find the Fourier transform of the invariant
distribution $\Lambda _\gamma$.  That is, find a linear functional
$\widehat{\Lambda _\gamma}$ such that
\[
\Lambda _\gamma (f)=\widehat{\Lambda _\gamma}\left( \hat{f}\right) ,
\]
where $\hat{f}$ is a function defined on the space of tempered
invariant ``eigendistributions'' on $G$. 
This space should include
the tempered irreducible characters of $G$ along with other invariant
distributions. 
 For example, if $\Pi$ is an
admissible representation of $G$ with character $\Theta_{\Pi}$, then
$$\hat{f}(\Pi) = \mathrm{tr}(\Pi(f)) = \int_G f(x)\Theta_{\Pi}(x)dx.$$
 The nature of the other distributions is an
intriguing problem.  The hope is that the Plancherel Formula for $G$
can be obtained through some limiting process for $\Lambda _\gamma$.

For example, if $G=SU(1,1)\cong SL(2,\R )$, we let
\[
\gamma =\begin{bmatrix}e^{i\theta _0} & 0 \\ 0 & e^{-i\theta _0}\end{bmatrix},\theta _0\neq 0,\pi .
\]
Then, $\gamma$ is a regular element in $G$, and
\[
G_\gamma =\T =\left\{\left. \begin{bmatrix}e^{i\theta} & 0 \\ 0 &
  e^{-i\theta}\end{bmatrix}\right|0\leq \theta <2\pi \right\} .
\]
After a simple computation, we get
\begin{align*}
F_f^\T (\gamma ) & =\left| e^{i\theta _0}-e^{-i\theta _0}\right|
\Lambda _\gamma (f) \\ & =-\frac{1}{2}\left( \pi ^{(+,+)}(f)-\pi
^{(+,-)}(f)\right) -\sum _{n\neq 0}\sgn (n)\chi _{\omega
  (n)}(f)e^{-in\theta _0} \\ & +\frac{i}{4}\left[ \int _\R \pi
  ^{(+,\nu )}(f)\frac{\sinh\left( \nu (\theta _0-\pi /2)\right)}{\sinh
    (\nu \pi /2)}d\nu-\int _\R \pi ^{(-,\nu )}(f)\frac{\cosh \left(
    \nu (\theta _0-\pi /2)\right)}{\cosh (\nu \pi /2)}d\nu\right] .
\end{align*}
The parameter $n\neq 0$ indexes the discrete series and the parameter
$\nu$ indexes the principal series representations of $G$.  The terms
$\pi ^{(+,+)}(f)$ and $\pi ^{(+,-)}(f)$ represent the characters of
the irreducible components of the reducible principal series, and we
obtain a ``singular invariant eigendistribution'' on $G$ by subtracting
one from the other and dividing by $2$.  This is exactly the invariant
distribution that makes harmonic analysis work.  It is called a
\emph{supertempered distribution} by Harish-Chandra.

This leads directly to the Plancherel Formula.  By a theorem of
Harish-Chandra, it follows that
\begin{align*}
\lim _{\theta \to 0}\left[ \frac{1}{i}\frac{d}{d\theta}\left[ F_f^{\T}(\gamma )\right] \right] & =8\pi f(1) \\
& =\sum _{n\in \Z}|n|\chi _{\omega (n)}(f)+1/2\int _0^\infty \pi ^{(+,\nu )}(f)\nu \coth (\pi /2\nu )d\nu \\
& +1/2\int _0^\infty \pi ^{(-,\nu )}(f)\nu \tanh (\pi /2\nu )d\nu .
\end{align*}

The representations of $SL(2,\R )$ were first determined by Bargmann
\cite{bargmann}.  In his 1952 paper \cite{hcsl2R}, Harish-Chandra gave
hints to the entire picture for Fourier analysis on real groups.  He
constructed the unitary representations, computed their characters,
found the Fourier transform of orbital integrals, and deduced the
Plancherel Formula.  This was done in about four and one-half pages.

We mention again that the support of the Fourier transform of the
tempered invariant distribution $\Lambda _\gamma$ contains not
only the characters of the principal series and the discrete series,
but also the tempered invariant distribution $$\frac{1}{2}\left( \pi
^{(+,+)}-\pi ^{(+,-)}\right).$$ This singular invariant
eigendistribution (appropriately normalized) is equal to $1$ on the
elliptic set and $0$ off the elliptic set, thereby having no effect on
harmonic analysis of the principal series.

Through the 1950s, along with an intensive study of harmonic analysis
on semisimple Lie groups, Harish-Chandra analyzed invariant
distributions, their Fourier transforms, and limit formulas related to
these.  This was mainly with reference to distributions on $C_c^\infty
(G)$.  He showed that $G$ has discrete series iff $G$ has a compact
Cartan subgroup.  For the rest of this section, we will assume that
$G$ has discrete series.  He also suspected quite early that the
irreducible unitary representations that occurred in the Plancherel
Formula would be indexed by a series of representations parameterized
by characters of conjugacy classes of Cartan subgroups.

In the 1960s, Harish-Chandra proved deep results about the character
theory of semisimple Lie groups, in particular, the discrete series
characters.  In developing the Fourier analysis on a semisimple Lie
group, he had to work with the smooth matrix coefficients of the
discrete series.  These matrix coefficients vanish rapidly at
infinity, but are not compactly supported.  This led to the definition
of the Schwartz space $\mathcal{C}(G)$ \cite{hcdiscrete2}.  The
Schwartz space was designed to include matrix coefficients of the
discrete series and slightly more.  The Schwartz space is dense in
$L^2(G)$, but is not contained in $L^1(G)$.  Moreover, the Schwartz
space $\mathcal{C}(G)$ does not contain the smooth matrix coefficients
of parabolically induced representations.  Nonetheless, the matrix
coefficients of these parabolically induced representations are
tempered distributions, that is, if $m$ is such a matrix coefficient
and $f\in \mathcal{C}(G)$, then $\int _Gfm$ converges.  Hence,
one can consider the orthogonal complement of these matrix
coefficients in $\mathcal{C}(G)$.

The collection of parabolically induced representations is indexed by
non-compact Cartan subgroups of $G$.  If $H$ is a Cartan subgroup of
$G$ with split component $A$, then the centralizer $L$ of $A$ is a
Levi subgroup of $G$.  Now the representations corresponding to $H$
are induced from parabolic subgroups with Levi component $L$, and the
subspace $\mathcal{C}_H(G)$ is generated by so called wave packets
associated to these induced representations.  Thus, we have an
orthogonal decomposition $$\mathcal{C}(G)=\bigoplus
\mathcal{C}_H(G),$$ where $H$ runs over conjugacy classes of Cartan
subgroups.  When $H$ is the compact Cartan subgroup of $G$,
$\mathcal{C}_H(G)$ is the space of cusp forms in $\mathcal{C}(G)$.
This decomposition of the Schwartz space is a version of the
Plancherel Theorem for $G$, and it is in this form that the Plancherel
Theorem appears for reductive $p$-adic groups.

As he approached his final version of the Plancherel Theorem and
Formula for real semisimple Lie groups, Harish-Chandra presented a
development of the Plancherel Formula for functions in $C_c^\infty
(G)$ in his paper \textit{Two Theorems on Semisimple Lie Groups}
\cite{hc2theorems}.  Here, he shows exactly how irreducible tempered
characters decompose the $\delta$ distribution.  In particular, for
$G$ of real rank $1$, he gives an explicit formula for the Fourier
transform of an elliptic orbital integral, and derives the Plancherel
Formula from this.  To understand the Plancherel Theorem for real
groups in complete detail, one should consult the three papers 
\cite{hcrrg1}, \cite{hcrrg2}, \cite{hcrrg3}, 
 and the expository renditions of this material
\cite{1970c}, \cite{1970a},  \cite{1972}.

\section{The Fourier Transform of Orbital Integrals, the Plancherel Formula, and Supertempered Distributions}

In a paper in \textit{Acta Mathematica} in 1973 \cite{sallywarner},
Sally and Warner re-derived, by somewhat different methods, the
inversion formula that Harish-Chandra proved in his ``Two Theorems''
paper \cite{hc2theorems}.  The purpose of the Sally--Warner paper was
to explore the support of the Fourier transform of an elliptic orbital
integral.  To quote: ``In this paper, we give explicit formulas for
the Fourier transform of $\Lambda _y$, that is, we determine a linear
functional $\widehat{\Lambda _y}$ such that
\[
\Lambda _y(f)=\widehat{\Lambda _y}\left( \hat{f}\right) ,f\in C_c^\infty (G).
\]
Here, we regard $\hat{f}$ as being defined on the space of
tempered invariant eigendistributions on $G$.  This space contains the
characters of the principal series and the discrete series for $G$
along with some `singular' invariant eigendistributions whose
character-theoretic nature has not yet been completely determined.''

In fact, the character theoretic nature of these singular invariant
eigendistributions was determined in a paper \cite{sallyherb} by Herb
and Sally in 1977.  In this paper, the present authors used results of
Hirai \cite{hirai}, Knapp--Zuckerman \cite{knappzuckerman}, Schmid
\cite{schmid}, and Zuckerman \cite{zuckerman} to show that, as in the
case of $SU(1,1)$, these distributions are alternating sums of
characters of limits of discrete series representations which can be
embedded as the irreducible components of certain reducible principal
series.  In his final published paper \cite{hcsupertemp}, Harish-Chandra
developed a comprehensive version of these singular invariant
eigendistributions, and he called them ``supertempered
distributions.''  These supertempered distributions include the
characters of discrete series along with some finite linear
combinations of irreducible tempered elliptic characters that arise
from components of reducible generalized principal series.  This
situation has already been illustrated for $SL(2,\R )$ in Section 2 of
this paper.  One notable fact about supertempered distributions is
that they appear discretely in the Fourier transforms of elliptic
orbital integrals; hence they play an essential role in the study of
invariant harmonic analysis.  For the remainder of this section, we
present a collection of results of the first author related to Fourier
inversion and the Plancherel Theorem for real groups.

In order to explain the steps needed to derive the Fourier transform
for orbital integrals in general, we first look in more detail at the
case that $G$ has real rank one.  In this case $G$ has at most two
non-conjugate Cartan subgroups: a non-compact Cartan subgroup $H$ with
vector part of dimension one, and possibly a compact Cartan subgroup
$T$.  We assume for simplicity that $G$ is acceptable, that is, the
half-sum of positive roots (denoted $\rho$) exponentiates to give a
well defined character on $T$.  The characters $\Theta ^T_{\tau }$ of
the discrete series representations are indexed by $\tau \in \hat T'$,
the set of regular characters of $T$, and the characters $\Theta
^H_{\chi }$ of the principal series are indexed by characters $\chi
\in \hat H$.  In addition, for $f \in C^{\infty}_c(G)$ we have
invariant integrals $F^T_f(t), t \in T$, and $F^H_f(a), a \in H$.
These are normalized versions of the orbital integrals $\Lambda
_{\gamma }(f), \gamma \in G$, which have better properties as
functions on the Cartan subgroups.

The analysis on the non-compact Cartan subgroup is elementary.  First,
as functions on $G'$, the set of regular elements of $G$, the
principal series characters are supported on conjugates of $H$.  In
addition, for $\chi \in \hat H, a \in H'=H \cap G'$, $\Theta ^H_{\chi
}(a)$ is given by a simple formula in terms of $\chi (a)$.  As a
result it is easy to show that the abelian Fourier transform $\hat
F^H_f(\chi ), \chi \in \hat H$, is equal up to a constant to $\Theta
^H_{\chi }(f)$, the principal series character evaluated at $f$.
Finally, $F^H_f \in C^{\infty}_c(H)$, and so the abelian Fourier
inversion formula on $A$ yields an expansion
\begin{equation} F^H_f(a) =  c_H \int _{\hat H}\epsilon (\chi )  \overline {\chi (a)}\Theta  ^H_{\chi }(f) d\chi , a \in H, \end{equation}
where $c_H$ is a constant depending on normalizations of measures and
$\epsilon (\chi ) = \pm 1$.

The situation on the compact Cartan subgroup is more complicated.
There are three main differences.  First, for $\tau \in \hat T', t \in
T' = T \cap G'$, $\Theta ^T_{\tau} (t)$ is given by a simple formula
in terms of the character $\tau (t)$.  However, $\Theta ^T_{\tau} $ is
also non-zero on $H'$.  Thus for $\tau \in \hat T', f \in
C_c^{\infty}(G)$, the abelian Fourier coefficient $\hat F^T_f(\tau )$
is equal up to a constant to $ \Theta ^T_{\tau} (f)$ plus an error
term which is an integral over $H$ of $F^H_f$ times the numerator of
$\Theta ^T_{\tau}$.  Second, the singular characters $\tau _0 \in \hat
T$ do not correspond to discrete series characters.  They do however
parameterize singular invariant eigendistributions $\Theta ^T_{\tau
  _0} $, and $\hat F^T_f(\tau _0 )$ can be given in terms of $\Theta
^T_{\tau _0}(f)$.  Finally, $F^T_f$ is smooth on $T'$, but has jump
discontinuities at singular elements.  Because of this there are
convergence issues when the abelian Fourier inversion formula is used
to expand $F^T_f$ in terms of its Fourier coefficients.

Sally and Warner were able to compute the explicit Fourier transform
of $F^T_f$ in the rank one situation where discrete series character
formulas on the non-compact Cartan subgroup were known.  The resulting
formula is very similar to the one for the special case of $SU(1,1)$
given in the previous section.  The discrete series characters and
singular invariant eigendistributions occur discretely in a sum over
$\hat T$ and the principal series characters occur in an integral over
$\hat A$ with hyperbolic sine and cosine factors.  They were also able
to differentiate the resulting formula to obtain the Plancherel
Formula.

The key to computing an explicit Fourier transform for orbital
integrals in the general case is an understanding of discrete series
character formulas on non-compact Cartan subgroups.  Thus we briefly
review some of these formulas.  The results are valid for any
connected reductive Lie group, but we assume for simplicity of
notation that $G$ is acceptable.  A detailed expository account of all
results about discrete series characters presented in this section is
given in \cite{herb2struct}.

Assume that $G$ has discrete series representations, and hence a
compact Cartan subgroup $T$, and identify the character group of $T$
with a lattice $L \subset E = i \mathfrak{t} ^*$.  For each $\lambda
\in E$, let $W(\lambda ) = \{ w \in W: w\lambda = \lambda \}$ where
$W$ is the full complex Weyl group, and let $E' = \{ \lambda \in E:
W(\lambda ) = \{1\} \}$.  Then $\lambda \in L' = L \cap E'$ is
regular, and corresponds to a discrete series character $\Theta
^T_{\lambda }$.  For $t \in T'$, we have the simple character formula
\begin{equation}  \Theta ^T_{\lambda }(t) =  \epsilon (E^+ ) \Delta (t)^{-1} \sum _{w \in W_K} \det (w) e^{w\lambda }(t), \end{equation}
where $\Delta $ is the Weyl denominator, $W_K$ is the subgroup of $W$
generated by reflections in the compact roots, and $ \epsilon (E^+ ) =
\pm 1$ depends only on the connected component (Weyl chamber) $E^+$ of
$E'$ containing $\lambda $.

Now assume that $H$ is a non-compact Cartan subgroup of $G$, and let $H^+$ be a connected component of $H'$.    Then for $h \in H^+$, 
\begin{equation}\label{dsc}  \Theta ^T_{\lambda }(h) =  c(H^+) \epsilon (E^+ ) \Delta (h)^{-1} \sum _{w \in W} \det (w)c(w\colon E^+\colon H^+)  \xi _{w, \lambda }(h), \end{equation}
where $c(H^+)$ is an explicit constant given as a quotient of certain
Weyl groups and the $c(w\colon E^+\colon H^+)$ are integer constants depending
only on the data shown in the notation.  The sum is over the full
complex Weyl group $W$, and for $w$ such that $c(w\colon E^+\colon H^+)$ is
potentially non-zero, $\xi _{w, \lambda }$ is a character of $H$
obtained from $w$ and $\lambda $ using a Cayley transform.  This
formula is a restatement of results of Harish-Chandra in
\cite{hcdiscrete1}.  In that paper, Harish-Chandra gave properties of
the constants $c(w\colon E^+\colon H^+)$ which characterize them completely.
These properties can in theory be used to determine the constants by
induction on the dimension of the vector component of $H$.  This
easily yields formulas when this dimension is one or two, but quickly
becomes cumbersome for higher dimensions.

With the above notation, it is easy to describe the singular invariant
eigendistributions corresponding to $\lambda \in L^s = L \backslash
L'$.  Let $\lambda_0 \in L^s$, and let $E^+$ be a chamber with
$\lambda _0 \in Cl( E^+)$.  The exponential terms $ \xi _{w, \lambda
  _0 }(h), h \in H^+$, still make sense, and the ``limit of discrete
series'' $\Theta ^T_{\lambda _0, E^+} = \lim _{\lambda \rightarrow
  \lambda_0, \lambda \in L \cap E^+} \Theta ^T_{\lambda }$ is given by
\eqref{dsc} using the constants from $E^+$.  Zuckerman
\cite{zuckerman} showed that the limits of discrete series are the
characters of tempered unitary representations of $G$.  The singular
invariant eigendistribution corresponding to $\lambda _0$ is the
alternating sum of the limits of discrete series taken over all
chambers with closures containing $\lambda _0$.
\begin{equation}  \Theta ^T_{\lambda _0} = [W(\lambda _0)]^{-1} \sum _{w \in W(\lambda _0)} \det w \  \Theta ^T_{\lambda _0, w E^+}. \end{equation}
     
 The main results of \cite{herbsallybig} are as follows.  Let $\Phi
 (\lambda _0)$ denote the roots of $T$ which are orthogonal to
 $\lambda _0$.  Then $\Theta ^T_{\lambda _0}$ vanishes if $\Phi
 (\lambda _0)$ contains any compact roots.  Thus we may as well assume
 that all roots in $\Phi (\lambda _0)$ are non-compact.  By using
 Cayley transforms with respect to the roots of $\Phi (\lambda _0)$ we
 obtain a Cartan subgroup $H$ and corresponding cuspidal Levi subgroup
 $M$.  Because the Cayley transform of $\lambda _0$ is regular with
 respect to the roots of $H$ in $M$, it determines a discrete series
 character of $M$, which can then be parabolically induced to obtain a
 unitary principal series character $\Theta ^H_{\lambda_0}$ of $G$.

\begin{thm} [Herb--Sally]    $ \Theta ^H_{\lambda _0} = \sum _{w \in W(\lambda _0 )} \Theta ^T_{\lambda_0, wE^+}$. 
 \end{thm}

It follows from Knapp \cite{knapp} that $\Theta ^H_{\lambda_0}$ has at
most $[W(\lambda _0)]$ irreducible components.  Thus each limit of
discrete series character is irreducible, and $\Theta ^T_{\lambda_0}$
is the alternating sum of the characters of the irreducible
constituents of $\Theta ^H_{\lambda_0}$.

In \cite{herbfourier}, Herb used the methods of Sally
and Warner, and the discrete series character formulas of
Harish-Chandra, to obtain a Fourier inversion formula for orbital
integrals for groups of arbitrary real rank.  As in the rank one case,
for any Cartan subgroup $H$ of $G$ we have normalized orbital
integrals $F^H_f(h), h \in H, f \in C^{\infty}_c(G)$.  We also have
characters $\Theta ^H_{\chi}, \chi \in \hat H$.  If $H$ is compact,
these are discrete series characters for regular $\chi $ and singular
invariant eigendistributions for singular $\chi $.  If $H$ is
non-compact, corresponding to the Levi subgroup $M$, then they are
parabolically induced from discrete series or singular invariant
eigendistributions on $M$.  Using standard character formulas for
parabolic induction, these characters can also be written using
Harish-Chandra's discrete series formulas for $M$.

Fix a Cartan subgroup $H_0$.  The goal is to find a formula
\begin{equation}\label{fif}   F^{H_0}_f(h_0) = \sum _H \int _{\hat H} \Theta ^H_{\chi}(f) K^H(h_0,\chi ) d\chi, \ h_0 \in H_0',  \end{equation}
 where $H$ runs over a set of representatives of conjugacy classes of Cartan subgroups of $G$,  $d\chi $ is Haar measure on $\hat H$, and $K^H(h_0 , \chi )$ is a function depending on $h_0$ and $\chi $.    
The problem is to compute the functions $K^H(h_0 , \chi )$, or at least show they exist.

 As in the rank one case, for $\chi _0 \in \hat H_0, f \in C_c^{\infty}(G)$, the abelian Fourier coefficient 
$\hat F^{H_0}_f(\chi _0)$ is equal up to a constant to $ \Theta ^{H_0}_{\chi _0} (f)$ plus an error term for each of the other Cartan subgroups.   The error term corresponding to $H$ is an integral over $H$ of the numerator of $\Theta ^{H_0}_{\chi _0}$ times $F^H_f$.    Because $\Theta ^{H_0}_{\chi_0}$ is parabolically induced, its character is non-zero only on Cartan subgroups of $G$ which are conjugate to Cartan subgroups of $M_0$, the corresponding Levi subgroup.     Thus the error term will be identically zero unless $H$ can be conjugated into $M_0$, but is not conjugate to $H_0$.     This implies in particular that the vector dimension of $H$ is strictly greater than that of $H_0$.     Thus if $H_0$ is maximally split in $G$ there are no error terms.    However, if $H_0=T$ is compact, then $M_0 =G$ and all non-compact Cartan subgroups contribute error terms.
  
Let $H$ be a Cartan subgroup of $M_0$ which is not conjugate to $H_0$ and let $M$ be the corresponding Levi subgroup.   In analyzing the error term corresponding to $H$, we obtain a primary term involving the characters $\Theta ^H_{\chi}(f), \chi \in \hat H$, plus secondary error terms, one for each Cartan subgroup of $M$ not conjugate to $H$.    This leads to messy bookkeeping, but the process eventually terminates since the vector dimension of the Cartan subgroups with non-zero error terms increases strictly at each step.    

In particular, if $H$ is a Cartan subgroup of $G$ not conjugate to a Cartan subgroup of $M_0$, then it never occurs in a non-zero error term and $K^H$ is identically zero.    Our original Cartan subgroup $H_0$ also is not involved in any error term, and we have
 \begin{equation}K^{H_0}(h_0, \chi_0) = c_{H_0}\epsilon (\chi _0) \overline{\chi _0(h_0)}, h_0 \in H_0', \chi _0 \in \hat H_0.\end{equation}
 The formulas for $K^H$ become progressively more complicated as the
 vector dimension of $H$ increases.  In particular, if $H$ is
 maximally split in $G$, then $K^H$ has contributions from error terms
 at many different steps.

Aside from the proliferation of error terms, the analysis which will
lead to the functions $K^H(h_0,\chi )$ involves two main problems that
do not occur in real rank one.  The main problem is that the final
formulas contain the unknown integer constants $c(w\colon E^+\colon H^+)$
appearing in discrete series character formulas.  These occur in
complicated expressions which can be interpreted as Fourier series in
several variables.  These series are not absolutely convergent and
have no obvious closed form.  Thus although \cite{herbsallybig} showed the
existence of the functions $K^H(h_0 , \chi )$, it does not result in a
formula which is suitable for applications.  In particular, it cannot
be differentiated to obtain the Plancherel Formula for $G$.  Second,
in the rank one case the analysis can be carried out for any $h \in
H'$.  However there are cases in higher rank, for example the real
symplectic group of real rank three, in which certain integrals
diverge for some elements $h \in H'$.  However, the analysis is valid
on a dense open subset of $H'$.
 
In order to improve these results and obtain a satisfactory Fourier
inversion formula similar to that of Sally and Warner for rank one
groups, it was necessary to have more information about the discrete
series constants.  The first of these improvements came from a
consideration of stable discrete series characters and stable orbital
integrals.

Assume that $G$ has a compact Cartan subgroup $T$, and use the notation
from the earlier discussion of discrete series characters.      For $\lambda \in L$ we define 
\begin{equation} \Theta ^{T,st}_{ \lambda } = [W_K]^{-1} \sum _{w \in  W} \Theta ^T_{w\lambda }.\end{equation}
If $\lambda \in L'$, then $\Theta ^{T,st}_{ \lambda }$ is called a stable discrete series character.    For $\lambda \in L^s$, we have
$\Theta ^{T,st}_{ \lambda } =0$.  
Similarly we define the stable orbital integral
\begin{equation} \Lambda ^{st}_t(f)  =  \sum _{w \in  W} \Lambda _{wt} (f), f \in C^{\infty}_c(G), t \in T'. \ \end{equation}
If we normalize the orbital integral as usual, we have
\begin{equation} F^{T, st}_f(t) = \Delta (t) \Lambda ^{st}_t(f)  =  \sum _{w \in  W} \det (w) F^T_f(wt).\end{equation}
Similarly, for any Cartan subgroup $H$ with corresponding Levi subgroup $M$ there is a series of stable characters  $\Theta ^{H,st}_{\chi }, \chi \in \hat H$, induced from stable discrete series characters of $M$.     We also obtain stable orbital integrals by averaging over the complex Weyl group of $H$ in $M$.    

Recall that there is a differential operator $\Pi$ such that
\begin{equation} f(1) = \lim _{t \rightarrow 1, t \in T'} \Pi F^T_f(t).  \end{equation}
Since the differential operator $\Pi$ transforms by the sign character
of $W$, it follows immediately that we also have
\begin{equation} f(1) = [W]^{-1} \lim _{t \rightarrow 1, t \in T'} \Pi F^{T,st}_f(t).  \end{equation}

The advantage of stabilizing is that the formulas for the stable
discrete series characters on the non-compact Cartan subgroups are
simpler than those of the individual discrete series characters.  The
Fourier inversion formula for stable orbital integrals involves only
these stable characters and has the general form
  \begin{equation}\label{sfif}   F^{T, st}_f(t) = \sum _H \int _{\hat H} \Theta ^{H, st}_{\chi}(f) K^{H,st}(t,\chi ) d\chi, \ t \in T'.  \end{equation}

When $G$ has real rank one the Fourier inversion formulas for the
stable orbital integrals are no simpler than those obtained by Sally
and Warner.  However when $G$ has real rank two there is already
significant simplification, and Sally's student Chao \cite{chao} was
able to obtain expressions for the functions $K^{H,st}(t,\chi )$ in
closed form and differentiate them to obtain the Plancherel Formula.

Herb \cite{herbmarseille}, \cite{herbamj} then developed the theory of
two-structures and showed that the constants occurring in stable
discrete series character formulas for any group can be expressed in
terms of stable discrete constants for the group $SL(2,\R)$ and the
rank two symplectic group $Sp(4,\R)$.  As a consequence she was able
to write each function $K^{H,st}(t,\chi )$ occurring in \eqref{sfif}
as a product of factors which occur in the corresponding formulas for
$SL(2,\R)$ and $Sp(4,\R)$.

This formula can be differentiated to yield the Plancherel Formula.
However, the Fourier inversion formulas for stable orbital integrals
are of independent interest, and much of the complexity of these
distributions is lost when they are differentiated and evaluated at
$t=1$.  In particular the functions occurring in the Plancherel
Formula, which had already been obtained by different methods by
Harish-Chandra \cite{hcrrg3}, reduce to a product of rank one factors
which occur in the Plancherel Formula for $SL(2,\R)$.  The discrete
series character formulas and Fourier inversion formula for $F^{T,
  st}_f(t)$ require both $SL(2,\R)$ and $Sp(4,\R)$ type factors coming
from the theory of two-structures.

In \cite{herbtams1} Herb was able to use Shelstad's ideas on endoscopy
to obtain explicit Fourier inversion formulas for the individual (not
stabilized) orbital integrals.  The idea is that certain weighted sums
of orbital integrals, $\Lambda ^{\kappa} _{\gamma }(f)$, correspond to
stable orbital integrals on endoscopic groups.  Thus their Fourier
inversion formulas can be computed as in \cite{herbamj}.  This is done
for sufficiently many weights $\kappa $ that the original orbital
integrals $\Lambda _{\gamma }(f)$ can be recovered.  Again, the theory
of two-structures was important, and the functions $K^H(h_0,\chi)$
occurring in \eqref{fif} can be given in closed form using products of
terms coming from the groups $SL(2,\R)$ and $Sp(4,\R)$.

Although this gave a satisfactory Fourier inversion formula, the
derivation is complicated by the use of stability and endoscopy.
Stability and endoscopy also combined to yield explicit, but
cumbersome, formulas for the discrete series constants $c(w\colon
E^+\colon H^+)$ occurring in \eqref{dsc}.  In \cite{herbtams2}, Herb
found simpler formulas for these constants that bypass the theories of
stability and endoscopy, and are easier to prove independently of
these results.  Using special two-structures called two-structures of
non-compact type, she obtained a formula for the constants $c(w\colon
E^+\colon H^+)$ directly in terms of constants occurring in discrete
series character formulas for $SL(2,\R)$ and $Sp(4,\R)$.  These
formulas could be used to give a direct and simpler proof of the
Fourier inversion formulas for orbital integrals given in
\cite{herbtams1}.

\section{The $p$-adic Case}

We now focus on the representation theory and harmonic analysis of
reductive $p$-adic groups.  Since the 1960s, there has been a flurry
of activity related to these groups.  Some of this has been generated
by the so-called ``Langlands Program'' (see Jacquet--Langlands
\cite{jacquetlanglands} and Langlands \cite{langlands}).  However, a
number of results in representation theory and harmonic analysis were
completed well before this activity related to the Langlands Program
by Bruhat \cite{bruhat}, Satake \cite{sataki}, Gel$'$fand--Graev
\cite{gelfandgraev}, and Macdonald \cite{macdonald}.  Of particular
interest were the results of Mautner \cite{mautnersph} that gave the
first construction of supercuspidal representations.  Here, a
supercuspidal representation is an infinite-dimensional, irreducible,
unitary representation with compactly supported matrix coefficients
(mod the center).  In the mid-1960s, for a $p$-adic field $F$ with
odd residual characteristic, all supercuspidal representations for
$SL(2,F)$ were constructed by Shalika \cite{shalikathesis}, and for
$PGL(2,F)$ by Silberger \cite{silberger}.  These two were Mautner's
Ph.D. students.  At roughly the same time, Shintani \cite{shintanisq}
constructed some supercuspidal representations for the group of
$n\times n$ matrices over $F$ whose determinant is a unit in the ring
of integers of $F$.  Shintani also proved the existence of a
Frobenius-type formula for computing supercuspidal characters as
induced characters.  Incidentally, in 1967--1968, the name
``supercuspidal'' had not emerged, and these representations were
called ``absolutely cuspidal,'' ``compactly supported discrete
series,'' and other illustrative titles.

We also note that, in this same period, Sally and Shalika computed the
characters of the discrete series of $SL(2,F)$ as induced characters
\cite{sallyshalika} (see also \cite{adss}), derived the Plancherel
Formula for $SL(2,F)$ \cite{sallyshalikaplanch}, and developed an explicit
Fourier transform for elliptic orbital integrals in $SL(2,F)$ 
\cite{sallyshalikaslnm}.  This Fourier transform led directly to the Plancherel
Formula through the use of the Shalika germ expansion \cite{shalikagerms}.
The guide for this progression of results was the 1952 paper of
Harish-Chandra on $SL(2,\R )$ \cite{hcsl2R}.

In the autumn of 1969, Harish-Chandra presented his first complete set
of notes on reductive $p$-adic groups \cite{vandijk}.  These are known
as the ``van Dijk Notes''.  These notes appear to be the origin of the
terms ``supercusp form'' and ``supercuspidal representation''.  They
present a wealth of information about supercusp forms, discrete series
characters, and other related topics.  At the end of the introduction,
Harish-Chandra states the following: ``Of course the main goal here is
the Plancherel Formula.  However, I hope that a correct understanding
of this question would lead us in a natural way to the discrete series
for $G$.  (This is exactly what happens in the real case.  But the
$p$-adic case seems to be much more difficult here.)''  It seems that
that Harish-Chandra favored the prefix ``super'' as in ``supercusp
form,'' ``supertempered distribution,'' etc.

We now proceed to the description of Harish-Chandra's Plancherel
Theorem (see \cite{hccollected4}) and Waldspurger's exposition of
Harish-Chandra's ideas \cite{waldsplanch}.  We then give an outline of
the current state of the discrete series of reductive $p$-adic groups
and their characters.  Finally, we give details (as currently known)
of the Plancherel Formula and the Fourier transform of orbital
integrals.

The background for Harish-Chandra's Plancherel Theorem was developed
in his Williamstown lectures \cite{hcwilliams}.  He showed that, using
the philosophy of cusp forms, one could prove a formula similar to
that for real groups that we outlined in Section 2.  He was able to do
this despite the lack of information about the discrete series and
their characters.

Following the model of real groups, for each special torus $A$,
Harish-Chandra constructed a subspace $\mathcal{C}_A(G)$ from the
matrix coefficients of representations corresponding to $A$.  These
representations are parabolically induced from relative discrete
series representations of $M$, the centralizer of $A$.  There are two
notable differences between the real case and the $p$-adic case.
First of all, because, in the $p$-adic case, there are discrete series
that are not supercuspidal (for example, the Steinberg representation
of $SL(2,F)$), the theory of the constant term must be modified.
Second, because of a compactness condition on the dual of $A$, it is
not necessary to consider the asymptotics of the Plancherel measure
that are required in the real case because of non-compactness.

Thus, even though the understanding of the discrete series and their
characters for $p$-adic groups is quite rudimentary, Harish-Chandra
succeeded in proving a version of the Plancherel Theorem.  This
version, as stated by Howe \cite{hccollected}, is: ``The (Schwartz)
space $\mathcal{C}(G)$ is the orthogonal direct sum of wave packets
formed from series of representations induced unitarily from discrete
series of (the Levi components of) parabolic subgroups $P$.  Moreover
if two such series of induced representations yield the same subspace
of $\mathcal{C}(G)$, then the parabolics from which they are induced
are associate, and the representations of the Levi components are
conjugate.''  Equivalently, as stated by Harish-Chandra (Lemma 5 of
\emph{The Plancherel Formula for Reductive $p$-adic Groups} in
\cite{hccollected4}), if $G$ is a connected reductive $p$-adic group
and $\mathcal{C}(G)$ is the Schwartz space of $G$,
then $$\mathcal{C}(G) = \sum_{A\in S} \mathcal{C}_A(G)$$ where $S$ is
the set of conjugacy classes of special tori in $G$ and the sum is
orthogonal.

In 2002, Waldspurger produced a carefully designed version of
Harish-Chandra's Plancherel Theorem.  This work is executed with
remarkable precision, and we quote here from Waldspurger's
introduction (the translation here is that of the authors of the
present article).

\ \\

``The Plancherel formula is an essential tool of invariant harmonic
analysis on real or $p$-adic reductive groups.  Harish-Chandra
dedicated several articles to it.  He first treated the case of real
groups, his last article on this subject being \cite{hcrrg3}.  A
little later, he proved the formula in the $p$-adic case.  But he
published only a summary of these results \cite{hccollected4}.  The
complete proof was to be found in a hand-written manuscript that was
hardly publishable in that state.  Several years ago, L.\ Clozel and
the present author conceived of a project to publish these notes.
This project was not realized, but the preparatory work done on that
occasion has now become the text that follows.  It is a redaction of
Harish-Chandra's proof, based on the unpublished manuscript.

\dots

As this article is appearing more than fifteen years after
Harish-Chandra's manuscript, we had the choice between scrupulously
respecting the original or introducing several modifications taking
account of the evolution of the subject in the meantime.  We have
chosen the latter option.  As this choice is debatable and the fashion
in which we observe the subject to have evolved is rather subjective,
let us attempt to explain the modifications that we have wrought.

There are several changes of notation: we have used those which seemed
to us to be the most common and which have been used since Arthur's
work on the trace formula.  We work on a base field of any
characteristic, positive characteristic causing only the slightest
disturbance.  We have eliminated the notion of the Eisenstein integral
in favor of the equivalent and more popular coefficient of the induced
representation.  We have used the algebraic methods introduced by
Bernstein.  They allow us to demonstrate more naturally that certain
functions are polynomial or rational, where Harish-Chandra proved
their holomorphy or meromorphy.  At the end of the article, we have
slightly modified the method of extending the results obtained for
semi-simple groups to reductive groups, in particular, the manner in
which one treats the center.  In fact, the principal change concerns
the `constant terms' and the intertwining operators.  Harish-Chandra
began with the study of the `constant terms' of the coefficients of
the induced representations and deduced from this study the properties
of the intertwining operators.  These latter having seemed to us more
popular than the `constant terms,' we have inverted the order, first
studying the intertwining operators, in particular their rational
extension, and having deduced from this the properties of the
`constant terms.'  All of these modifications remain, nevertheless,
minor and concern above all the preliminaries.  The proof of the
Plancherel formula itself (sections VI, VII and VIII below) has not
been altered and is exactly that of Harish-Chandra.''

\ \\

It remains to address the current status of the three central problems of
harmonic analysis on reductive $p$-adic groups.  These are the
construction of the discrete series, the determination of the
characters of the discrete series, and the derivation of the Fourier
transform of orbital integrals as linear functionals on the space of
supertempered distributions.

  There is a long list of authors who have attacked the construction
  of discrete series of $p$-adic groups over the past forty years.  We
  limit ourselves to a few of the major stepping stones.  The work of
  Howe \cite{howetame} on $GL(n)$ in the tame case set the stage for a
  great deal of the future work.  Howe's supercuspidal representations
  for $GL(n)$ were proved to be exhaustive by Moy in \cite{moy}.
  Further work in the direction of tame supercuspidals may be found in
  the papers  \cite{morris1} and \cite{morris2} of L. Morris.
  
In the mid 1980s, Bushnell and Kutzko attacked $GL(n)$ in the wild
case.  Their main weapon was the theory of types, and the definitive
results for $GL(n)$ and $SL(n)$ were published in
\cite{bushnellkutzko}, \cite{bksl1}, and \cite{bksl2}.  While in the
tame case, one gets a reasonable parameterization in terms of
characters of tori, it does not seem that such a parameterization can
be expected in the wild case.  It is difficult to associate certain
characters with any particular torus, as well as to tell when
representations constructed from different tori are distinct.  We  also
mention the work of Corwin on division algebras in both the tame
\cite{corwintame} and the wild \cite{corwinwild} case.

A big breakthrough came in J.-K.~Yu's construction of tame
supercuspidal representations for a wide class of groups in
\cite{jkyu}.  In this paper, Yu points to the fact that he was guided
by the results of Adler \cite{adlerrefined} at the beginning of this
undertaking.  Under certain restrictions on $p$, Yu's supercuspidal
representations were proved to be exhaustive by Ju-Lee Kim
\cite{jlkim} using tools from harmonic analysis in a remarkable way.
Throughout this period, the work of Moy--Prasad \cite{moyprasad1},
\cite{moyprasad2} was quite influential. Also, Stevens
\cite{stevens} succeeded in applying the Bushnell--Kutzko methods to
the classical groups to obtain all their supercuspidal representations as
induced representations when the underlying field has odd residual
characteristic.  Finally, major results have been obtained by M{\oe}glin
and Tadic for non-supercuspidal discrete series in
\cite{moeglintadic}.  There is still much work to be done, but
considerable progress has been made.

The theory of characters has been slower in its development.  There
are two avenues of approach that have been cultivated.  The first is
the local character expansion of Harish-Chandra.  If $\mathcal{O}$ is
a $G$-orbit in $\mathfrak{g}$, then $\mathcal{O}$ carries a
$G$-invariant measure denoted by $\mu_\mathcal{O}$ (see, for example,
\cite{rao}).  The Fourier transform of the distribution $f\mapsto
\mu_\mathcal{O}(f)$ is represented by a function
$\widehat{\mu_\mathcal{O}}$ on $\mathfrak{g}$ that is locally summable
on the set of regular elements $\mathfrak{g}'$ in $\mathfrak{g}$.  The
local character expansion is:
\begin{thm}
  Let $\pi$ be an irreducible smooth representation of $G$.  There are
  complex numbers $c_{\mathcal{O}}(\pi)$, indexed by nilpotent orbits
  $\mathcal{O}$, such that $$\Theta_\pi(\exp Y) = \sum_{\mathcal{O}}
  c_\mathcal{O}(\pi)\widehat{\mu_{\mathcal{O}}}(Y)$$ for $Y$
  sufficiently near $0$ in $\mathfrak{g}'$.
\end{thm}
This result is presented in Harish-Chandra's Queen's Notes
\cite{hcqueens} and is fully explicated in \cite{debackersally}.  The
local character expansion could be a very valuable tool if three
problems are overcome.  These are: (1) determine the functions
$\widehat{\mu_\mathcal{O}}$, (2) find the constants
$c_{\mathcal{O}}(\pi)$, and (3) determine the domain of validity of the
expansion.  For progress in these directions, see Murnaghan
\cite{murnaghan1}, \cite{murnaghan2}, Waldspurger \cite{walds2001},
 DeBacker--Sally \cite{debackersallygerms}, and DeBacker \cite{debackerhomo}.

The second approach is the direct use of the Frobenius formula for
induced characters to produce full character formulas on the regular
elements in $G$.  See Harish-Chandra \cite{vandijk} (p.\ 94), Sally
\cite{sallyosaka}, and Rader--Silberger \cite{radersilberger}.  This
approach has been used by DeBacker for $GL(\ell)$, $\ell$ a prime
\cite{debackerthesis}, and Spice for $SL(\ell)$, $\ell$ a prime
\cite{spicesll}.  Recent work of Adler and Spice \cite{adlerspice} and
DeBacker and Reeder \cite{debackerreeder} shows some promise in this
direction, but their results are still quite limited.  The
paper \cite{adlerspice} of Adler and Spice gives an interesting report
on the development and current status of character theory on reductive
$p$-adic groups.  For additional results on the theory of characters,
consult the papers of Cunningham and Gordon \cite{cunningham} and
Kutzko and Pantoja \cite{kutzkopantoja}.

We finish this paper with an update on the Plancherel Theorem, the
Plancherel Formula, and the Fourier transform of orbital integrals in
the $p$-adic case.  As regards the Plancherel Theorem, it seems that
some flesh is beginning to appear on the bones.  Thus, for some
special cases, an explicit Plancherel measure related to the
components in the Schwartz space decomposition has been found (see
Shahidi \cite{shahidi1}, \cite{shahidi2}, Kutzko--Morris
\cite{kutzkomorris}, and Aubert--Plymen \cite{aubertplymen1},
\cite{aubertplymen2}). The results seem to be applicable mainly to
$GL(n)$ and $SL(n)$.  In some cases, restrictions on the residual
characteristic have been completely avoided.  These methods seem to
a great extent to be independent of explicit character formulas.  It
would be interesting to determine how far these techniques can be
carried for general reductive $p$-adic groups.

It is one of the purposes of this paper to point out the nature of the
Plancherel Formula in the theory of harmonic analysis on reductive
$p$-adic groups.  As was the case originally with Harish-Chandra, the
Plancherel Formula should be considered as the Fourier transform of
the $\delta$ distribution regarded as an invariant distribution on a
space of smooth functions on the underlying group.  This is achieved
in the real case by determining the Fourier transform of an elliptic
orbital integral and applying a limit formula involving differential
operators to deduce an expression for $f(1)$ as a linear functional on
the space of tempered invariant distributions.  This space is directly
connected to the space of tempered irreducible characters of $G$ along
with some additional supertempered virtual characters.  It appears to
be the case that, to accomplish this goal, one has to have a full
understanding of the irreducible tempered characters of $G$.  This, of
course, requires a detailed knowledge of the discrete series.  This is
exactly the approach that was detailed in Section 3.

As pointed out by Harish-Chandra, a complete knowledge of the discrete
series and their characters would yield the Plancherel measure for
$p$-adic groups exactly as in the real case.  In the $p$-adic case, the
role of differential operators in the limit formula to obtain $f(1)$
is assumed by the \emph{Shalika germ expansion}.

%\subsection{Shalika Germs}

\ \\
\noindent {\bf Shalika Germs}

For a connected semi-simple $p$-adic group $G$, Shalika defines in
\cite{shalikagerms} $$I_f(x) = \int_{G(x)} fd\mu,$$ where $x$ is a
regular element in $G$, $G(x)$ is its conjugacy class, $\mu$ is a
$G$-invariant measure on $G(x)$, and $f\in C_c^\infty(G)$.  Shalika
shows that $I_f(x)$ has an asymptotic expansion in terms of the
integrals $$\Lambda_{\mathcal{O}}(f) = \int_{\mathcal{O}} fd\mu$$ of
$f$ over the unipotent conjugacy classes $\mathcal{O}$.  Here, for
$\mathcal{O}=\{1\}$, we take $\Lambda_{\mathcal{O}}(f)=f(1)$.  The
coefficients $C_{\mathcal{O}}(x)$ occurring in this expansion are
called the \emph{Shalika germs}.

We start with $G=SL(2,F)$ where $F$ has odd residual characteristic,
and then use Shalika germs to produce the Plancherel Formula for $G$.
This result of Sally and Shalika was proved in 1969 and is presented
in detail in \cite{sallyshalikaslnm}.  We repeat it here to indicate
the role that such a formula can play in the harmonic analysis on a
reductive $p$-adic group.

Let $T$ be a compact Cartan subgroup of $G$.  For each nontrivial
unipotent orbit $\mathcal{O}$, there is a subset $T_\mathcal{O}$ of
the set of regular elements in $T$ such that the following asymptotic
expansion holds.
$$F_f^T(t) =|D(t)|^{1/2}I_f(t)\sim -A_T|D(t)|^{1/2}f(1) + B_T
\sum_{\dim{\mathcal{O}}>0} C_\mathcal{O}(t)\Lambda_{\mathcal{O}}(f)$$
where the Shalika germ $C_\mathcal{O}(t)$ is the characteristic
function of $T_\mathcal{O}$.  The constants $A_T$ and $B_T$ depend on
normalization of measures and whether $T$ is ramified or unramified.

  By summing products of characters, we are led to the following
  expression.
  \begin{eqnarray*}
    \mu(T)I_f(t) &=& \sum_{\Pi\in
      D}\overline{\chi_{\Pi}(t)}\hat{f}(\Pi) +
    \frac{1}{2}\sum_{\Pi\in\mathrm{RPS}_V}\overline{\chi_\Pi(t)}\hat{f}(\Pi)
    \\ &-& \frac{q+1}{2q}\mu(A_1) \int_{\!\!\!\!\tiny\begin{array}{c}
      \xi\in\widehat{F^\times} \\ \xi|A_{h_0+1}=1\end{array}}
    |\Gamma(\xi)|^{-2} \hat{f}(\xi)d\xi \\ &+&
    \frac{q}{2}\mu(A_1)\kappa_{T}|D(t)|^{-1/2}
    \int_{\!\!\!\!\tiny\begin{array}{c} \xi\in\widehat{F^\times}
      \\ \xi|A_{h_0+1}=1\end{array}} \hat{f}(\xi)d\xi
  \end{eqnarray*}
  This is the Fourier transform of the elliptic orbital integral
  corresponding to the regular element $t$.  Note the occurrence of
  the characters of the reducible principal series, denoted
  $\mathrm{RPS}_V$, corresponding to the three sgn characters on
  $F^\times$.  As in the case of $SL(2,\R)$, each represents the
  difference of two characters divided by $2$, and that difference is
  $0$ except on the compact Cartan subgroups corresponding to the sgn
  character associated to the quadratic extension $V$.  So again,
  these singular tempered invariant distributions (see
  \cite{ramseysally}) appear in the Fourier transform of an elliptic
  orbit.

  Using Shalika germs, we are led directly to the Plancherel Formula for
  $SL(2,F)$.
  $$\mu(K)f(1) = \sum_{\Pi\in D}\hat{f}(\Pi)d(\Pi) +
  \frac{1}{2}\left(\frac{q^2-1}{q}\right)\mu(A_1)\int_{\xi\in
    \widehat{F^\times}}|\Gamma(\xi)|^{-2}\hat{f}(\xi)d\xi$$
  
It is clear that a complete theory of the Fourier transform of orbital
integrals would lead to direct results about lifting, matching, and
transferring orbital integrals.  Such a theory would entail a deep
understanding of discrete series characters and their properties.  A
start in this direction may be found in papers of Arthur
\cite{arthuretc}, \cite{arthurft} and Herb \cite{herbell},
\cite{herbsuper}.  We expect to return to this subject in the near
future.

%\bibliographystyle{plain}
%\bibliography{msnbib}

\def\cprime{$'$} \def\cprime{$'$}

\end{document}